\documentclass[reqno]{amsart}

\usepackage{amsmath}
\usepackage{amsfonts}
\usepackage{rotating}
\usepackage{bbm}
\usepackage[all,cmtip]{xy}
\usepackage{amscd}
 
\usepackage{tikz-cd}              
\usepackage{tikz}
 \usepackage{pgfplots}

 \usetikzlibrary{arrows}
 \tikzset{>=latex}
\usetikzlibrary{arrows.meta}
\pgfplotsset{compat=1.10}
\usetikzlibrary{decorations.markings}

\usepgfplotslibrary{fillbetween}
\usetikzlibrary{patterns}
\usepackage{mathrsfs}
\usepackage{color}
\usepackage{xcolor}
\usepackage{mathtools}
\usepackage{amssymb}
\usepackage[latin1]{inputenc}
\usepackage{graphicx}
\usepackage{dsfont}
\usepackage{color}
\usepackage{hyperref}
\usepackage{enumerate}
\usepackage{marginnote}

\usepackage{subcaption}
 \usepackage{comment}

\newtheorem{theorem}{Theorem}[section]

\newtheorem{lemma}[theorem]{Lemma}
\newtheorem{proposition}[theorem]{Proposition}

\theoremstyle{definition}

\theoremstyle{remark}
\newtheorem{remark}[theorem]{Remark}

\newcommand{\p}{\partial}

\newcommand{\C}{{\mathbb{C}}}
\newcommand{\I}{{\mathbb{I}}}
\newcommand{\R}{{\mathbb{R}}}

\newcommand{\N}{{\mathbb{N}}}
\newcommand{\HH}{{\mathbb{H}}}

\renewcommand{\epsilon}{\varepsilon}
\renewcommand{\theta}{\vartheta}

\usepackage{graphicx}
 
\begin{document}

\title[]{  Spatial Stark-Zeeman systems and their regularizations}

\author[S. Kim]{Seongchan Kim}
   \address{ Department of Mathematics Education, Kongju National University, Gongju 32588, Republic of Korea}
   \email {seongchankim@kongju.ac.kr}

 \author[K.Ruck ]{Kevin Ruck}
   \address{ Seoul National University, Department of Mathematical Sciences, Research institute in Mathematics, Gwanak-Gu, 
			Seoul 08826, South Korea }
   \email {kevin.ruck@snu.ac.kr}

\setcounter{tocdepth}{3}

\date{Recent modification; \today}
 \begin{abstract}
  
In this article, we study spatial Stark-Zeeman systems which describe the dynamics of a charged particle moving in three-dimensional space under the influence of a Coulomb potential, a magnetic field, and an electric field, possibly time-dependent. Such systems are modeled by  Hamiltonian flows on the cotangent bundle of an open subset of $\R^3, $ equipped with a twisted symplectic structure.  The presence of the Coulomb singularity leads to the study of collision orbits, and hence understanding the regularization of these orbits is essential for global dynamical properties.  
 We   investigate   regularization techniques for  spatial Stark-Zeeman systems, both in time-independent and time-dependent cases. In particular, in the time-dependent case, following a new regularization method developed by Barutello, Ortega, and Verzini, we   formulate the  corresponding regularized variational principles and  carefully analyze the effects of the magnetic and electric terms under the Kustaanheimo-Stiefel transformation.  The resulting regularized action functional      yields a variational characterization of collision orbits and facilitates further analysis of periodic solutions. 
  Our results provide  a   general scheme for regularizing spatial Stark-Zeeman systems, opening the door for further applications in symplectic geometry, Floer theory, and celestial mechanics.

  \end{abstract}
   \maketitle

   


\section{Introduction}

The Stark-Zeeman systems provide a rich class of Hamiltonian systems describing the motion of a charged particle under the influence of a Coulomb force, a magnetic field, and an electric field. The interaction between the Coulomb singularity, the magnetic field that   twists the symplectic structure, and the presence of  an electric field  causes significant analytical challenges. In particular, the appearance of collision orbits should be very carefully regularized to  be able to  analyze global dynamical properties and the existence of periodic orbits.

The planar Stark-Zeeman systems, i.e.~those, whose configuration space is an open subset in the plane,  have been studied, for example, in \cite{CFvK17Jplus, CFZ23SZ, Fra25BOV, FraKimtwocenterSZ}.
By contrast, the spatial case, where the configuration space is the three-dimensional space, exhibits new dynamical and geometric complexity.  In this article, we develop a framework for the regularization of spatial Stark-Zeeman systems, in particular, using the Kustaanheimo-Stiefel transformation. In the case where the electric field is time-dependent, following the recent regularization technique due to Barutello, Ortega, and Verzini \cite{BOV21}, we obtain a regularized variational principle such that the critical points correspond to collisional periodic orbits in the original problem. This may lay  the groundwork for further analysis of periodic orbits using variational and Floer-theoretic techniques.


Assume that the simply connected open set $\mathfrak{U}_0 \subset \R^3$ contains the origin and let
\[
\mathfrak{U} = \mathfrak{U}_0 \setminus \{ 0 \}.
\]
For a smooth function $E\colon S^1 \times \mathfrak{U}_0 \to \R$ we write
\[
E_t\colon \mathfrak{U}_0 \to \R, \quad E_t(q)=E(t,q), \quad t \in S^1.
\]
The function $E$  may be constant in time.  We now define the potential
\[
V \colon S^1 \times \mathfrak{U} \to \R, \quad q \mapsto V (t,q):= - \frac{1}{\lvert q \rvert}+ E (t,q) 
\]
and the Hamiltonian
\begin{equation*}\label{eq:Hamiltoniantwisted}
H \colon S^1 \times T^*\mathfrak{U} \to \R, \quad H(t, q,p):= \frac{1}{2}\lvert p \rvert^2 + V(t,q).
\end{equation*}
As before,   we write $V_t(q) = V(t,q) $ and $H_t(q,p)=H(t,q,p).$ We now consider a closed two-form
\[
\sigma_{B} = \frac{1}{2}\sum_{i,j=1}^3 B_{ij}(q) dq_i \wedge dq_j \in \Omega^2(\mathfrak{U}_0), \quad {\rm where} \;\; B_{ij} = -B_{ji},
\]
referred to as a {\it magnetic form}, and  define the twisted symplectic form by
\[
\omega_\sigma := \sum_{j=1}^3 dp_j \wedge dq_j + \pi^*\sigma_B   \in \Omega^2( T^*\mathfrak{U}_0),
\]
where  $\pi \colon T^*\mathfrak{U}_0 \to \mathfrak{U}_0 $ denotes the footpoint projection.  Note that we follow the usual notation in symplectic geometry by denoting the base point coordinates by $q_j$ and the fiber coordinates by $p_j$. 

The twisted Hamiltonian vector field $X_{H,\sigma}$ is implicitly given by
\[
\omega_\sigma( \cdot , X_{H,\sigma}) = dH_t. 
\]
A periodic solution  $x \colon S^1 \to T^*\mathfrak{U}$ then solves the first-order ODE
\[
\dot{x}(t) = X_{H,\sigma}(x(t)), \quad t \in S^1.
\]
 Since the twisted Hamiltonian vector field is explicitly written as
  \[
 X_{H,\sigma}(x(t)) = \sum_{i=1}^3 \left[ p_i(t) \partial_{q_i}  + \left(\sum_{j=1}^3 B_{ij}(q(t))p_j(t) - \partial_{q_i}V_t(q(t))    \right)\partial_{p_j}      \right]
 \]
   for a periodic solution $x=(q,p)\colon S^1 \to T^*\mathfrak{U}$ the loop $q \colon S^1 \to \mathfrak{U}$ satisfies the second-order ODE
 \begin{equation}\label{eq:secondorder}
 \ddot{q} = B(q)\dot{q} - \nabla V_t(q) = B(q)\dot{q} - \frac{q}{\lvert q \rvert^3} - \nabla E_t(q).
 \end{equation}

Alternatively, instead of  describing the magnetic field via a twisted symplectic structure one can incorporate it  {directly} into the Hamiltonian itself by what in physics is known as \textit{minimal coupling}: From Maxwell's equations we know that the magnetic field needs to satisfy $d\sigma_B=0,$ and since $\mathfrak{U}_0$ is simply connected, this implies that there exists a one-form $A\in \Omega^1(\mathfrak{U}_0)$ with 
\[
\sigma_B=dA.
\]
In local coordinates  {this} implies
\begin{equation}\label{eq:relationbetweenAandB}
B_{ij} = \frac{\partial A_j}{\partial q_i } - \frac{\partial A_i}{\partial q_j}.
\end{equation}
Using this one-form we  define a fiber translation on $T^*\mathfrak{U}_0$ via
\[
\text{T}_A \colon  T^*\mathfrak{U}_0 \to T^*\mathfrak{U}_0, \quad \alpha_q\mapsto \alpha_q-A(q),
\]
 called the minimal coupling for the magnetic field $B,$ and the Hamiltonian 
\[
H_A:=(\text{T}_A)^*H =\frac{1}{2}\lvert p - A(q)\rvert^2 + V(t,q),
\]
called the \textit{minimally coupled} Hamiltonian. Note that  $\text{T}_A$ has an obvious inverse  $\text{T}_{-A}$. Pulling back the standard symplectic structure $\omega_0$ on $T^*\mathfrak{U}_0$ via $\text{T}_{-A}$ gives
\[
(\text{T}_{-A})^*\omega_0 = \omega_\sigma,
\]
{which implies that
\[
    \text{T}_{-A} \colon  \left(T^*\mathbb{R}^3,\omega_\sigma,H\right) \to \left(T^*\mathbb{R}^3,\omega_0,H_A\right)
\]
is a symplectomorphism between the two Hamiltonian systems.}
This shows that the two ways  {we introduced} of describing a magnetic field in Hamiltonian mechanics are equivalent. Further, it also shows that the specific one-form $A$ we use to define the minimal coupling does not matter to the physics of the system as long as it satisfies $\sigma_B=dA$. This freedom in the choice of $A$ is called \textit{gauge freedom}.

Let us denote by
\[
\left< \cdot, \cdot \right> \colon \mathfrak{L} \R^3 \times \mathfrak{L}\R^3 \to \R
\]
the $L^2$-inner product on the free loop space $\mathfrak{L}\R^3=C^\infty(S^1, \R^3),$ defined by
\[
\left< \xi, \eta\right>= \int_0^1 \left< \xi(t), \eta(t)\right>dt, \quad \xi, \eta \in \mathfrak{L}\R^3.
\]
The associated $L^2$-norm is indicated by 
\[
 \| \xi \| =\sqrt{ \left< \xi, \xi \right>} , \quad \xi \in  \mathfrak{L}\R^3  .
 \]
The Lagrangian action functional corresponding to the Hamiltonian $H_A$ is then given by 
\begin{equation}\label{eq:functionalgeneralSZ}
\mathscr{A} \colon \mathfrak{L}\mathfrak{U} \to \R, \quad \mathscr{A}(q):=\frac{1}{2}\| \dot{q}\|^2 + \int_{S^1}q^* A - \int_0^1 V_t(q(t))dt.
\end{equation}
A direct computation shows that    critical points of the functional $\mathscr{A}$ correspond to the solutions of the Newtonian equation \eqref{eq:secondorder}.  

In the study of the second-order ODE \eqref{eq:secondorder}, collision orbits naturally arise. To extend solutions smoothly through such collisions, it becomes necessary to regularize the~system.   In the  time-independent setting, this is a classical topic in celestial mechanics, with a variety of well-established local regularization techniques, see, for instance,   \cite{Birk15regularization, KS65regularization, LC20regularization}. In Section \ref{sec:timeindpe} we apply the Moser and Kustaanheimo-Stiefel maps to regularize time-independent   Stark-Zeeman systems. 
However, in the time-dependent case, the total energy is no longer conserved, and hence the classical approach fails.   This motivates the need for a new regularization approach. To address this, we develop a {\it Kustaanheimo-Stiefel type BOV-regularization} by adapting a method recently introduced  by Barutello-Ortega-Verzini \cite{BOV21}, which provides   a smooth  framework for handling collision orbits and enables a  variational analysis of periodic solutions to non-local differential equations.  A detailed account of  different types of  BOV-regularizations for planar Stark-Zeeman systems can be found in  \cite{Fra25BOV, FraKimtwocenterSZ}.

Following \cite[Definition 1.1]{BOV21}, we define:

\smallskip

\noindent
{\bf Definition.}
A continuous map   $q \colon S^1 \to \R^3$   is said to be a {\it generalized solution} of \eqref{eq:secondorder} if it satisfies the following conditions:
\begin{enumerate}
    \item[$i)$] the zero set $\mathcal{Z}_q = \{ t \in S^1 \mid q(t) =0\}$ is finite;
    
    \item[$ii)$]   the map $q$ is smooth and satisfies \eqref{eq:secondorder} on the complement of $  \mathcal{Z}_q$;  
    
    \item[$iii)$]  the energy 
    \[
  \mathscr{E}(t):=   \frac{1}{2}\lvert \dot{q}(t)\rvert ^2 -\frac{ 1  }{\lvert q(t) \rvert}   + \int_t^1 \dot{E}_s(q(s))ds + E_t(q(t)),  \quad t \in S^1 \setminus \mathcal{Z}_q
    \]
    admits a continuous extension on the whole  circle $S^1,$  so that   it is constant along every generalized solution.

\end{enumerate}
\noindent

\smallskip

The main result of this article is the following assertion which informally summarizes Theorem \ref{them:main11}. 

\smallskip

\noindent
{\bf Theorem A.} {\it A Kustaanheimo-Stiefel type BOV-regularization provides a variational approach to the analysis of generalized solutions in spatial Stark-Zeeman systems. }

\smallskip

\subsection*{Acknowledgments}
The second author was supported by the National Research Foundation of Korea under the grants RS-2023-00211186, NRF-2020R1A5A1016126 and (MSIT) RS-2023-NR076656.

\medskip

\section{Examples of spatial Stark-Zeeman systems}  Before diving deeper into the analysis of   spatial Stark-Zeeman systems and their regularizations, we first present several  important examples of   Stark-Zeeman systems  {that should} serve  as motivation. While the general Stark-Zeeman system is typically described using terminology such as electric   and magnetic fields, many of the most relevant examples arise in the context of celestial mechanics. Indeed, all of the examples  discussed below are systems studied   exclusively within celestial mechanics.

\subsection{ The rotating Kepler problem}
The most fundamental example of a spatial Stark-Zeeman system in celestial mechanics is the spatial rotating Kepler problem. Its Hamiltonian is given by
\begin{equation*}\label{H_rotKepler}
    H(q,p):= \frac{1}{2}|p|^2 - \frac{1}{\lvert q\rvert} +p_1q_2-p_2q_1, \quad (q,p) \in T^*(\R^3 \setminus \{0\}),
\end{equation*}
where the first two terms are just the usual Kepler problem and the remaining part is inducing a rotation of the frame of reference by an angular velocity of $1$. 
After rearranging this Hamiltonian as
\[
H (q,p ) = \frac{1}{2}\left((p_1+q_2)^2+(p_2-q_1)^2+p_3^2 \right) - \frac{1}{\lvert q\rvert} -\frac{1}{2}\lvert q\rvert^2
\]
we can use equation~\eqref{eq:functionalgeneralSZ} to write the Lagrangian action functional as
\begin{equation*}
    \begin{aligned}
        \mathscr{A}_{\text{RKP}}(q)=&\int_0^1 \frac{1}{2}\lvert \dot{q}(t)\rvert^2-\dot{q}_1(t)q_2(t)+\dot{q}_2(t)q_1(t)+\frac{1}{\lvert q(t)\rvert}+\frac{1}{2}\lvert q(t)\rvert ^2 dt\\
        =& \frac{1}{2}\|\dot{q}\|^2+\int_{S^1}q^*\left(-x_2dx_1+x_1dx_2\right) - \int_0^1 \left(-\frac{1}{\lvert q(t)\rvert}-\frac{1}{2}\lvert q(t)\rvert^2\right) dt.
    \end{aligned}
\end{equation*}
From this we can easily see that the rotating Kepler problem is a Stark-Zeeman system as described in the introduction. The one-form $A$ associated with the magnetic field is here given by
\[
A(x)= -x_2dx_1+x_1dx_2
\]
and describes the Coriolis force that is induced by the rotating frame. The electric potential is in this case constant in time and is given by
\[
E(x)= -\frac{1}{2}|x|^2.
\]
This term is due to the centrifugal force experienced by a particle in this system. 

\medskip

The remaining examples we will present in this section all build upon this idea: the Coulomb potential is the Kepler potential, the magnetic term describes the Coriolis force induced by a rotating frame and the electric potential is a combination of centrifugal force  and the gravitational attraction of other celestial bodies.

\subsection{The circular restricted three-body problem} The (spatial) circular restricted three-body problem (CR3BP)   describes the motion of  a massless particle moving under the gravitational influence of two massive bodies that move around their shared center of gravity in a circular motion according to the solution of the two-body problem (i.e.~ignoring the influence of the particle). If we denote by $a(t)$ and $b(t)$ the time-dependent positions of the two bodies, we can write the Hamiltonian of the CR3BP as
\[
H(q,p,t)= \frac{1}{2}\lvert p\rvert^2-\frac{\mu}{\lvert q-a(t)\rvert}-\frac{1-\mu}{\lvert q-b(t)\rvert},
\]
where $\mu$ is the mass of body $a$ and we normalized the combined mass of body $a$ and $b$ to be $1$. Assuming that   the barycenter lies at the origin, this implies that $|a(t)|=1-\mu$ and $|b(t)|=\mu$ for all $t$. In this form the Hamiltonian does not really fit into the format of a Stark-Zeeman system, since both Kepler potentials are time-dependent. However, when we again consider a rotating frame around the center of gravity of the two bodies, the potential becomes time-independent:
\[
H^\prime(q,p)= \frac{1}{2}\lvert p\rvert^2 +p_1q_2-p_2q_1-\frac{\mu}{\lvert q-a\rvert}-\frac{1-\mu}{\lvert q-b\rvert }
\]
Note that we can always choose the coordinate system in such a way that the movement of the two bodies lies in the $x$-$y$ plane.
One problem that remains is that both Kepler potentials have   singularities, hence we cannot simply view them as the electric potential $E$.  First we have to restrict our system to energies below the first critical energy value, where the components of the Hill region around each body $a$ and $b$ are disjoint and bounded. Now choose one of the bodies: we pick $a$ for this discussion. Then we restrict the position space of the system to $\mathfrak{U}$, which is the component of the Hill region around body $a$ at an energy below the first critical energy. This allows us to view the Kepler potential of $b$ as well-defined function $E\colon \mathfrak{U}_0\to \mathbb{R}$. After a coordinate shift that puts body $a$ into the center, we have the following Hamiltonian:
\[
H_{\text{CR3BP}}(q,p)= \frac{1}{2}|p|^2 +p_1q_2-p_2\left(q_1+(1-\mu)\right)-\frac{\mu}{|q|}-\frac{1-\mu}{|q+\vec{e}_1|}
\]
Here, $\vec{e}_1$ stands for the vector $(1,0,0)$. Note that we can assume that the positions $a$ and $b$ lie on the $x$-axis with $a_1=1-\mu$ and $b_1=-\mu$. Now we can rearrange the Hamiltonian as in the previous section and get the Lagrangian action functional in accordance with equation~\eqref{eq:functionalgeneralSZ}:

\begin{align*}
\mathscr{A}_{\text{CR3BP}}(q)= & \; \frac{1}{2}||\dot{q}||^2 +\int_{S^1}q^*\left(-x_2dx_1+\left(x_1+(1-\mu)\right)dx_2\right) \\
& \; - \int_0^1 \left(-\frac{\mu}{|q|}-\frac{1}{2}\left((q_1+(1-\mu))^2 + q_2^2\right)-\frac{1-\mu}{|q+\vec{e}_1|}\right) dt.
\end{align*}
This shows that the CR3BP is also a (time-independent) Stark-Zeeman system with the data
\[
A(x)= -x_2dx_1+\left(x_1+(1-\mu)\right)dx_2
\]
and
\[
E(q)= -\frac{1}{2}\left((q_1+(1-\mu))^2 + q_2^2\right)-\frac{1-\mu}{|q+\vec{e}_1|}.
\]
 \begin{remark}
     Note that following the original definition of a Stark-Zeeman system, we defined the potential $V$ as 
     \[
        V(t,q)=-\frac{1}{|q|}+E(t,q).
    \]
    However, in applications to celestial mechanics, the Coulomb term usually carries an additional constant denoting the mass ratio between the bodies involved. This additional constant will of course not spoil any properties that are valid for the original definition of a Stark-Zeeman system.
 \end{remark}

\subsection{The bicircular restricted four-body problem}
So far, our examples are only time-independent versions of Stark-Zeeman systems, which can be regularized using the already well-established standard local methods (see Section~\ref{sec:timeindpe}). The bicircular restricted four-body problem (BR4BP) is now the first example that will lead to a  time-dependent Stark-Zeeman system for which the  local regularizations are not applicable.

This system generalizes the CR3BP mentioned above. A typical example of the CR3BP is the Earth-Moon-spaceship configuration. However, in reality, the Earth-Moon system  itself orbits the Sun, and   especially close to the Moon, the gravitational influence of the Sun is not   negligible compared to that  of the Earth. To account for this, the BR4BP was introduced.

\begin{figure}[h]
    \centering
    \begin{tikzpicture}
    \draw[arrows=->] (-4,0) --(4,0);
    \draw[arrows=->] (0,-3) --(0,3);
    \draw (0,0) node {$\times$};
    \draw (0.2,-0.3) node {$B_1$};
    \draw (-2.5,0) node {$\bullet$};
    \draw (-2.5,-0.3) node {$a$};
    \draw (1.5,0) node {$\bullet$};
    \draw (1.5,-0.3) node {$b$};
    \draw[dashed] (0,0) -- (1.5,2.6);
    \draw[blue] (0, 0) to[out=90, in=270] (0.75, 0.3);
    \draw[blue] (1.5, 0) to[out=90, in=270] (0.75, 0.3);
    \draw[blue] (0.75, 0.4) node {$\mu$};
    \draw[blue] (0, 0) to[out=90, in=270] (-1.25, 0.3);
    \draw[blue] (-2.5, 0) to[out=90, in=270] (-1.25, 0.3);
    \draw[blue] (-1.25, 0.4) node {$1-\mu$};
    \draw[dashed] (0, 0) circle[radius=3];
    \draw (1.25,2.17) node {$\times$};
    \draw (1.6,2) node {$B_2$};
    \draw (1.5,2.6) node {$\bullet$};
    \draw (1.7,2.8) node {$s$};
    \draw[green!40!black] (1.5,2.6) to[out=150, in=330] (1.205, 2.48);
    \draw[green!40!black] (1.25,2.17) to[out=150, in=330] (1.205, 2.48);
    \draw[green!40!black] (0.75, 2.65) node {$l-\nu$};
    \draw[green!40!black] (0, 0) to[out=150, in=330] (0.281, 1.28);
    \draw[green!40!black] (1.25,2.17) to[out=150, in=330] (0.281, 1.28);
    \draw[green!40!black] (0.14, 1.4) node {$\nu$};
    \draw[red!60!black] (2,0) arc (0:60:2);
    \draw[red!60!black] (1.5,0.8) node {$\alpha$};
    \draw[arrows=->] (1.5,2.6) --(2.37,2.1);
    \draw (2.3,2.5) node {$\omega_r$}; 
\end{tikzpicture}
    \caption{The bicircular restricted four-body problem}
    \label{fig:BCR4BP}
\end{figure}
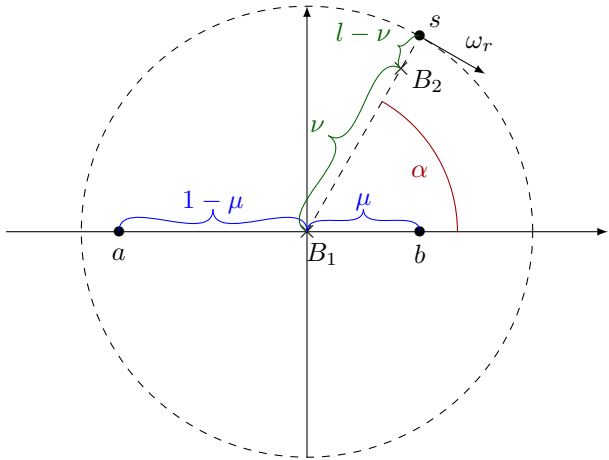
The setting is as follows: two bodies $a$ and $b$ are moving around their barycenter $B_1$ in  circular motion with  normalized angular velocity of $1$. In addition,   a third body $s$   moves together with   $B_1$ around their barycenter $B_2,$ also in   circular motion, with angular velocity $\omega_{B_2}.$ The  distances of $B_1$ and $s$ from $B_2$ are $\nu$ and $l - \nu,$ respectively,  where $l$ denotes  the relative distance from the $a$-$b$ system to $s$. For simplicity, we assume that all these motions lie in the same plane.   Within  this setup, we consider a spaceship  bound to body $a,$ which  moves under the gravitational influence of all three bodies $a$, $b$ and $s$. To express the system in the form of a  Stark-Zeeman system a sequence of coordinate changes is required. These transformations are not Galilean and thus introduce  a different kind of inertial forces.

First we shift the coordinates so that the barycenter $B_1$ lies at the origin. This causes the barycenter $B_2$ to rotate around the origin and induces  Coriolis and centrifugal forces depending additionally on $\nu$ and $\omega_{B_2}$. Then we introduce a rotating frame like for the CR3BP to fix the bodies $a$ and $b$ at their position on the $x$-axis. This induces the typical Coriolis and centrifugal forces seen in the previous sections and causes the body $s$ and the barycenter $B_2$ to move at the new relative angular velocity of $\omega_r$. In the last step we then shift the coordinate system such that $a$ lies at the origin.  Overall, this description leads to the Hamiltonian
\begin{equation*}
    \begin{aligned}
        \tilde{H}_{\text{BC}}(q,p)=& \frac{1}{2}|p|^2 +\omega_{B_2}p_1\left(q_2-\nu\sin(\omega_r t+\alpha)\right)-\omega_{B_2}p_2\left(q_1+a_1-\nu\cos(\omega_r t+\alpha)\right)\\ &+p_1q_2-p_2\left(q_1+a_1\right)-\frac{\mu}{|q|}-\frac{1-\mu}{|q+a-b|}-\frac{m_s}{|q+a-s(t)|},
    \end{aligned}
\end{equation*}
where $m_s$ is the relative mass of body $s$ with respect to the combined mass of $a$ and $b$, $\alpha$ is the angle determined by the starting configuration of $a$, $b$ and $s$ (see Figure~\ref{fig:BCR4BP}) and $s(t)$ the time-dependent position of the body $s$ given by
\[
s(t):= l\begin{pmatrix}\cos(\omega_r t+\alpha)\\\sin(\omega_r t+\alpha)\end{pmatrix}.
\]
To separate the Coriolis potential (magnetic term) from the centrifugal potential (part of the electric potential)  we define 
\[
\tilde{A}_t(x):= \begin{pmatrix}\omega_{B_2}\left(x_2-\nu\sin(\omega_r t+\alpha)\right)+x_2\\ -\omega_{B_2}\left(x_1+a_1-\nu\cos(\omega_r t+\alpha)\right)-x_1\\0
\end{pmatrix}, \quad E_\text{cf}(x)=\vec{e}_3\times \tilde{A}_t(x).
\]
Here $\vec{e}_3$ denotes the vector $(0,0,1)$. Since we are working on $\mathbb{R}^6$ we can view $\tilde{A}_t$ as a one-form and rewrite the above Hamiltonian as
\[
     \tilde{H}_{\text{BC}}(q,p,t)= \frac{1}{2}\left|p-\tilde{A}_t(q)\right|^2 -\left|E_\text{cf}(q)\right|^2 -\frac{\mu}{|q|}-\frac{1-\mu}{|q+a-b|}-\frac{m_s}{|q+a-s(t)|}.
\]
As we can see now, the Hamiltonian has a time-dependent magnetic term, which does not conform to our definition of a Stark-Zeeman system. However, we also know that the Coriolis force for a rotation in the plane only depends on the angular momentum and not on the position of the rotation axis. Hence, we should be able to use the \textit{gauge freedom} in the choice of the magnetic one-form to eliminate this `non-physical' time-dependence. First, let us calculate the magnetic two-form corresponding to $\tilde{A}_t$:
\begin{align*}
    \sigma_B=& \; d\tilde{A}_t(x)\\
    =& \;  d\left(\left(\omega_{B_2}\left(x_2-\nu\sin(\omega_r t+\alpha)\right)+x_2\right)dx_1\right.\\
    & \; \left. -\left(\omega_{B_2}\left(x_1+a_1-\nu\cos(\omega_r t+\alpha)\right)+x_1\right)dx_2\right)\\
    =&  \; \left(\omega_{B_2}+1\right)dx_2\wedge dx_1-\left(\omega_{B_2}+1\right)dx_1\wedge dx_2\\
    =&  \; -2\left(\omega_{B_2}+1\right)dx_1\wedge dx_2
\end{align*}
This calculation shows that indeed we can choose a different one-form $A$ that is time-independent. We will choose
\[
A(x)= \left(\omega_{B_2}+1\right)x_2 dx_1-\left(\omega_{B_2}+1\right)x_1 dx_2
\]
Now we can proceed as in the previous sections and use equation~(\ref{eq:functionalgeneralSZ}) to  deduce the relevant Lagrangian action functional
\begin{align*}
\mathscr{A}_{\text{BC}}(q)= &\frac{1}{2}||\dot{q}||^2 +\int_{0}^1 \left\langle A(q(t)),\dot{q}(t)\right\rangle dt -\int_0^1 -\frac{\mu}{|q|}dt\\ &- \int_0^1 \left(-\left|E_\text{cf}(q)\right|^2  -\frac{1-\mu}{|q+a-b(t)|}-\frac{m_s}{|q+a-s(t)|}\right) dt.
\end{align*}
Since $A$ is now time-independent, this functional is in accordance with the definition of a Stark-Zeeman system, where the time-dependent electric potential is given by
\[
E_t(x)= -\left|E_\text{cf}(x)\right|^2 -\frac{1-\mu}{|x+a-b(t)|}-\frac{m_s}{|x+a-s(t)|}.
\]
Note that we have to restrict the domain $\mathfrak{U}_0$ to be the Hill regions for energies below the first critical energy value, in order for the function $E_t$ to be well-defined and smooth everywhere. 

\begin{remark}
    As described at the beginning of this section the BR4BP is defined by making the CR3BP rotate around an additional body. Hence, we should be able to recover the CR3BP from the above Hamiltonian when we set the mass $m_s$ of the additional body and the angular velocity of the second barycenter $\omega_{B_2}$ to zero. In the above formula we can easily see that this will change the Hamiltonian of the BR4BP to $H_\text{CR3BP}$.
\end{remark}

 \subsection{Celestial mechanics beyond Stark-Zeeman}
 So far, all the   examples discussed above--which are variants of  restricted $N$-body problem--turned out to be   Stark-Zeeman systems. Hence, one might suspect that nearly all systems in celestial mechanics  fall into this category. However, in this section, we  present a brief non-example to illustrate that one does not need to stray far from the standard models to encounter  a system in celestial mechanics that does not conform to the framework of  a Stark-Zeeman system.

The setting is that of a restricted three-body problem, hence, in inertial coordinates $\left(q^\prime,p^\prime\right)$ the Hamiltonian reads
\[
H(q^\prime,p^\prime,t)= \frac{1}{2}\lvert p^\prime\rvert^2-\frac{\mu}{\lvert q^\prime-a(t)\rvert}-\frac{1-\mu}{\lvert q^\prime-b(t)\rvert}.
\]
But now each of the two bodies $a$ and $b$ move on a general ellipse (instead of a circle as before), where the two ellipses share the same focus  at the center of mass. The orbit of the bodies is then given by
\[
    a(t)=(1-\mu)r(t)\begin{pmatrix} \cos(f(t)\\ \sin(f(t))\\0
    \end{pmatrix}
\]
with true anomaly $f(t)$ and distance of $a$ from the focus 
\[
r(t)=\frac{1-e^2}{1+e\cos(f(t))}.
\]
 Similar to the CR3BP we can introduce a non-inertial frame in which the two primaries are fixed. For the elliptic case these are the pulsating coordinates (see \cite[Chapter~8.10]{Meyer17}). Further, it is also beneficial to change the time variable to be the true anomaly $f$. Overall the resulting Hamiltonian is given by
 \begin{align*}
     H_\text{ER3BP}(q,p,f)= & \; \frac{1}{2}|p+\vec{e}_3\times q|^2 -\frac{\mu\cdot r(f)}{\lvert q\rvert}-\frac{(1-\mu)\cdot r(f)}{\lvert q+\vec{e}_1\rvert}\\
     & \; -\left(1-\frac{r^\prime(f)}{r(f)}\right) \left|q+(1-\mu)\vec{e}_1\right|^2,
 \end{align*}
 called the elliptic restricted three-body problem (ER3BP).
 
 If we now take a look at the Coulomb potential, we see that it contains the time-dependent term $r(f)$. Since we introduced pulsating coordinates, which periodically scale the distance in the system, it should be not surprising that we also have to scale the gravitational attraction accordingly. But this means there should be no way to get  rid  of the time-dependence in the Coulomb potential without changing the physics of the system. Therefore, the ER3BP is not a Stark-Zeeman system.

\medskip

\section{Kustaanheimo-Stiefel transformation}

\subsection{Kustaanheimo-Stiefel map}\label{KSmap}

Following \cite{Zhao15KS}, we denote by
\[
\HH=\{z = z_0 + z_1 i + z_2 j + z_3 k \mid (z_0, z_1, z_2, z_3) \in \R^4 \} , \quad \I \HH = \{ z \in \HH \mid {\rm Re}(z) =0\} 
\]
skew-field of quaternions and the subset of purely imaginary quaternions, respectively. Here, the real and the imaginary parts of a quaternion of $z = z_0 + z_1 i + z_2 j + z_3 k \in \HH$ are respectively given by
\[
{\rm Re}(z) = z_0 \quad  \text{ and } \quad {\rm Im}(z) = z_1 i + z_2 j + z_3 k.
\]
Note that, in a straightforward manner, $\HH$ and $\I\HH$ are isomorphic to $\R^4$ and $\R^3,$  respectively.

The conjugation of $z \in \HH$ is defined to be
\[
\bar{z} = z_0 - z_1 i - z_2 j - z_3 k .
\]
We then define the Kustaanheimo-Stiefel   map (KS map) by
\begin{equation}\label{eq:KSmap1}
    \Phi\colon \HH \to \I \HH, \quad \Phi(z) = \bar{z} i z.
\end{equation}
This map is an $S^1$-covering map with respect to the $S^1$-action
\[
\Phi(e^{i\theta}z) = \Phi(z), \quad \theta \in S^1,\; z \in \HH.
\]
To see this consider the following lemma.
\begin{lemma}\label{S^1cover}
    The quotient of the KS map with respect to the above $S^1$-action 
    \[
        \hat{\Phi}\colon \HH/S^1 \to \I \HH, \qquad [z]\mapsto \Phi(z)
    \]
    is a bijection.
\end{lemma}
\begin{proof}
 We fix $q\in \I\HH$. Recall that   every three-dimensional rotation in $\I\HH\cong\mathbb{R}^3$ corresponds to   conjugation with a   unit quaternion. Thus, we can find a unit quaternion $u$ such that $uq\bar{u}=|q|i$. Now define $z:= \sqrt{|q|}$, which gives us 
    \[
        \bar{z}iz=|q|i=uq\bar{u}.
    \]
    Rearranging this we get
    \[
        q=\overline{zu}izu,
    \]
    which shows the surjectivity of $\hat{\Phi}.$

    Next, assume that there are two quaternions $a$ and $b$ such that 
    \[
    \bar{a}ia=\bar{b}ib.
    \]
    This equation is equivalent to
    \[
        i=\frac{\overline{b\bar{a}}}{|a|^2}i\frac{b\bar{a}}{|a|^2}.
    \]
    Since $a$ and $b$   have the same norm, $ {b\bar{a}}/{|a|^2}$ is a unit quaternion. Using the polar form of the quaternions it is now an  easy exercise to show that the only unit quaternions that satisfy $i=\bar{z}iz$ are of the form $e^{i\theta}$. Hence, we can conclude that $b=e^{i\theta}a$ for some $\theta \in \R$. This completes the proof.
\end{proof}

By lifting this map   we obtain  
\[
\mathfrak{P} \colon T^*(\HH \setminus \{0\}) \to \I \HH \times \HH, \quad (z, w) \mapsto \left(q=\Phi(z), \; p=\frac{\bar{z}iw}{2\lvert z\vert^2}  \right)
\]
Note that the map $\mathfrak{P}$ satisfies the relation  
\[
\mathfrak{P}^*(pdq)=\frac{1}{2}\sum\limits_{i=1}^{3}\left(w_idz_i-z_idw_i\right)
\]
and is invariant under the $S^1$-action as well:
\[
 \mathfrak{P} (e^{i \theta}z, e^{i \theta}w )=\mathfrak{P} (z,w), \quad   \theta \in S^1, \;(z,w) \in T^*(\HH \setminus \{0\}).
\]

\subsection{Kustaanheimo-Stiefel transformation}

Abbreviate by $\mathfrak{L}\HH=C^\infty(S^1, \HH)$ the free loop space of $\HH.$
We denote by
\[
\left< z_1, z_2 \right> = \int_0^1 \left< z_1(\tau), z_2(\tau) \right> d\tau = \int_0^1 {\rm Re}(\bar{z}_1(\tau)z_2(\tau))d\tau
\]
the $L^2$-inner product of $z_1, z_2 \in L^2(S^1, \HH).$ The associated $L^2$-norm of $z \in L^2(S^1, \HH)$ is defined as
\[
\| z\|:= \sqrt{ \left< z ,z \right>}.
\]

Suppose that $z \colon S^1 \to \HH$  is a continuous map  with finite zero set
\begin{equation*}\label{eq:conditionofz}
\mathcal{Z}_z:=z^{-1}(0).
\end{equation*}
We define a reparametrization map
\begin{equation}\label{eq:tz}
t_z \colon S^1 \to S^1, \quad t_z(\tau) = \frac{1}{\| z \|^2} \int_0^\tau \lvert  z(s)\rvert^2 ds.
\end{equation}
Since the zero set $\mathcal{Z}_z$ is finite,   the map $t_z$ is strictly increasing, from which we see that it is a homeomorphism of $S^1.$ Denote by $\tau_z:=t_z^{-1} \colon S^1 \to S^1$ its continuous inverse.  We now define $q_z \colon S^1 \to \I \HH$ by
\begin{equation}\label{eq:definitionofqz}
q_z(t):=\Phi( z(\tau_z(t))),
\end{equation}
where $\Phi$ denotes the KS map, see \eqref{eq:KSmap1}. Note that 
\begin{align*}
\int_0^1 \frac{1}{\lvert q_z(t) \rvert}dt  = \frac{1}{\| z \|^2}< \infty.
\end{align*}

Conversely, let $q \colon S^1 \to \I\HH$ be a continuous map with finite zero set 
\[
\mathcal{Z}_q :=q^{-1}(0) 
\]
satisfying
\[
\int_0^1\frac{1}{\lvert q(t)\rvert}dt < \infty.
\]
We define a reparametrization map
\begin{equation*}\label{eq:tauq}
\tau_q \colon S^1 \to S^1, \quad \tau_q(t):= \left( \int_0^1 \frac{1}{\lvert q(s)\rvert}ds\right)^{-1} \int_0^t\frac{1}{\lvert q(s)\rvert}ds
\end{equation*}
As before, since the zero set $\mathcal{Z}_q$ is finite, the map $\tau_q$ is a homeomorphism whose inverse $t_q:=\tau_q^{-1}$ has the derivative
\[
t_q'(\tau) = \left( \int_0^1 \frac{1}{\lvert q(s)\rvert}ds \right) q(t_q(\tau)), \quad \tau \in S^1.
\]
Define $z_q \colon S^1 \to \HH$ by
\begin{equation}\label{eq:fromqtozPhi}
\Phi( z_q(\tau)) = q(t_q(\tau)). 
\end{equation}
Since the zero set $\mathcal{Z}_{z_q}$ of $z_q$ satisfies the relation $\mathcal{Z}_{z_q} = \tau_q(\mathcal{Z}_q),$ it is in particular finite. Therefore, we can define the reparametrization map $t_{z_q}\colon S^1 \to S^1$ as in \eqref{eq:tz}. Denote by $\tau_{z_q}$ its inverse. Arguing as in \cite[Section 2.2]{CFV23-1} one can show that
\[
t_{z_q} = t_z \quad \quad \text{ and } \quad \quad \tau_{z_q} = \tau_q, 
\]
implying that $q$ is the Kustaanheimo-Stiefel transformation of $z_q$ given in \eqref{eq:definitionofqz}.

\section{Time-independent cases}\label{sec:timeindpe}

In this section we assume that the potential $V$ is time-independent and discuss the  regularization procedure for the corresponding Stark-Zeeman systems.

Since energy is preserved by the flow, we consider the  energy level
\[
\Sigma_c = H_A^{-1}(c), \quad c \in \R 
\]
and the associated Hill region
\[
\mathfrak{K}_c :=\pi(\Sigma_c)= \{ q \in \mathfrak{U} \mid V(q)\leq c\} .
\]
Due to the singular term $-1/\lvert q \rvert$ in    the potential $V,$   the Hill region $\mathfrak{K}_c$ contains a unique connected component whose closure contains the origin. We denote this distinguished component by  
\[
\mathfrak{K}_c^b \subset \mathfrak{K}_c.
\]
We define  
\[
c_0:=\sup\{ \tilde c \mid {\text{conditions (1) and (2) below hold for all $c<\tilde{c}$}}\}\in \R \cup \{+\infty\},
\]
where 
\begin{enumerate}
\item every $c<\tilde{c}$ is a regular value of $H_A$;
    \item  for every $c<\tilde{c},$ the component $\mathfrak{K}_c^b$ is diffeomorphic to a closed unit disk with the origin removed.
    
\end{enumerate}

\subsection{Moser regularization}

In the time-independent case the standard regularization procedure is based on finding a suitable compactification of the energy hypersurface $\Sigma_c$ for a chosen energy $c$, such that the dynamics on the compactification coincides with the dynamics on $\Sigma_c$ outside the collision locus.

In Moser regularization this compactification is achieved by mapping the   momenta, which escape to infinity as we approach the collision locus, onto the sphere via the inverse stereographic projection such that the north pole corresponds to the point at infinity. To be more precise, we first choose an energy level $c<c_0$ at which we would like to regularize. 
If we would  simply map the momenta in the spatial Stark-Zeeman system to the 3-sphere, it would not be a symplectic transformation and the dynamics on the hypersurface would not need to be preserved. Hence, we first apply the symplectic switch map
\[
\mathfrak{sw}:T^*\mathbb{R}^3 \to T^*\mathbb{R}^3, \qquad (q,p)\mapsto (-p,q).
\]
Now the base manifold is describing the momenta and the cotangent fiber the positions. This allows us to consider the cotangent lift of the inverse stereographic projection
\[
    \mathfrak{\phi}: \mathbb{R}^3 \to S^3 \subset \mathbb{R}^4, \quad (x_1,x_2,x_3)\mapsto \left(\frac{2x_1}{1+\lvert x\rvert^2}, \frac{2x_2}{1+ \lvert x\rvert^2}, \frac{2x_3}{1+\lvert x\rvert^2}, \frac{ \lvert x\rvert^2-1}{1+\lvert x\rvert^2}\right).
\]
If we denote 
\[
    \Phi:= T^*\phi \circ \mathfrak{sw}: T^*\mathbb{R}^3 \to T^*S^3
\]
the regularized Hamiltonian is given by
\[
    H_M:= \Phi^*\left(|q|(H-c)\right),
\]
where $H$ is the unregularized Hamiltonian, and the full Hamiltonian system is $\left(T^*S^3, \omega_{\text{can}}, H_M\right)$.  {Heuristically speaking,} this regularization is compactifying the energy hypersurface $\Sigma_c$ by adding an $S^1$ family of momenta at the collision locus, which represents the directions of the momenta at which the collisions happen. The regularized collision orbits are trajectories that `bounce' back from the collision locus with minus the momentum they have at the time of the collision. For more details on why this regularization leads to a well-defined compactified energy hypersurface for all classes of Stark-Zeeman systems we refer to \cite[Section 3.5]{CFvK17Jplus}.

\begin{remark} To apply techniques from symplectic and contact geometry to the study of the dynamics on  the regularized energy level, denoted by $\overline{\Sigma}_c$, it is of particular importance to determine whether   $\overline{\Sigma}_c$ is of contact type. It was shown in  \cite{CJK203BP}   that the Liouville vector field $\eta \partial_{\eta}$ (with $\eta$ denoting   the fiber coordinate)  is transverse to the regularized energy level $\overline{\Sigma}_c$ for  {all energies below the first critical energy}   in the CR3BP.

\end{remark}

\subsection{Kustaanheimo-Stiefel regularization}
Another way of regularizing a spatial time-independent Stark-Zeeman system is via the lift $\mathfrak{P}$  of the Kustaanheimo-Stiefel map (see Section~\ref{KSmap}). Since the image of this map is $\I\HH\times\HH,$ but $H_A$ is only defined on $\mathbb{R}^6 \cong \I\HH \times \I \HH  $, we first have to restrict $\mathfrak{P}$ to a suitable subset. For this we define the map
\[
    BL\colon \HH\times\HH \to \mathbb{R}, \qquad (z,w) \mapsto \text{Re}(\bar{z}iw).
\]
Note that the set $BL^{-1}(0)$ is exactly the subset in the domain of $\mathfrak{P}$ that is mapped to $\I\HH\times\I\HH=\mathbb{R}^6$. Hence, we consider the set
\[
   \Gamma_1:= BL^{-1}(0) \cap T^*\left(\HH\setminus\{0\}\right).
\]
For a given regular value $c$ of $H_A$ we then define the Kustaanheimo-Stiefel regularized Hamiltonian
$K_c \colon \Gamma_1 \to \R$ by
\begin{align*}
  K_c (z,w) &= \left(\mathfrak{P}\big\vert_{\Gamma_1}\right)^*\left(\lvert q \rvert ( H_A -c)\right)(z,w)\\
  &= \frac{1}{8}\frac{\lvert \bar{z}iw \rvert^2}{\lvert \bar{z}iz \rvert} +\lvert \bar{z}iz \rvert^2\left( E(\bar{z}i z)  - c +\frac{1}{2}\lvert  A(\bar{z}iz) \rvert^2   \right) - \frac{1}{2}\left< \bar{z}i w ,  \mathcal{A}(z)  \right> - 1.
\end{align*}
Note that we can extend this Hamiltonian to $T^*\HH$ by simply extending the norm and vector product to their four-dimensional equivalent:
\[
   K_c (z,w)  = \frac{1}{8}\lvert w \rvert^2 +\lvert z \rvert^2\left( E(\bar{z}i z)  - c +\frac{1}{2}\lvert  A(\bar{z}iz) \rvert^2   \right) - \frac{1}{2}\left< \bar{z}i w ,  \mathcal{A}(z)  \right> - 1 .
\]
It is easy to see that this Hamiltonian does not have any divergences and carries a globally well-defined dynamics. However, when it comes to relating the regularized dynamics to the unregularized one we have an obvious problem: the former has eight degrees of freedom and the latter only six. One of the additional degrees of freedom can be explained by the $S^1$-symmetry that $\mathfrak{P}$ artificially introduces into the regularized system. The other one is due to extending $K_c$ from $\Gamma_1$ to $T^*\HH$. But this means that only those Hamiltonian trajectories in $T^*\HH$ that in addition lie in the set $BL^{-1}(0)$ are related to trajectories in the unregularized system, the other ones are just a non-physical byproduct of the regularization procedure.  

To finish this section we want to briefly discuss the topology of the resulting regularized energy hypersurface. First consider  the case of the Kepler problem where we have $E\equiv 0$ and $A\equiv 0.$  Then the Hamiltonian $K_c$ reads
\begin{equation*} 
K_c^{\rm Kepler}(z,w) = \frac{1}{8}\lvert w \rvert^2   - c  \lvert z \rvert^2    - 1
\end{equation*}
To find the regularized energy level, we first note  that the zero level set $(K_c^{\rm Kepler})^{-1}(0)$ is diffeomorphic to the seven sphere $S^7.$
Its intersection  with the set $BL^{-1}(0)$ is $ S^3 \times S^3$ (see \cite{Zhao15KS} for further details).

Now go back to the general case: We can homotop our regularized energy hypersurface through hypersurfaces of lower and lower energy levels by decreasing $c$ until we reach a very negative energy. To reach an energy hypersurface of the Kepler problem we then continuously turn off the external electric and magnetic fields. This homotopy is well-defined as long as we do not cross a critical energy of the corresponding Hamiltonian. Hence, it follows that for every energy below the first critical energy value the bounded component of the KS regularized energy level of a spatial Stark-Zeeman system is diffeomorphic to $S^3 \times S^3$. Note that, different from \cite{Zhao15KS}, we do not divide the energy hypersurface by the $S^1$-action.

\begin{remark}
    Similar to the Levi-Civita regularization (see \cite[Remark 3.7]{CFvK17Jplus}) we can find an $S^1$-covering map from the Kustaanheimo-Stiefel regularized energy hypersurface to the one given by the spatial Moser regularization.
\end{remark}

\medskip

\section{Time-dependent cases}

We now  study time-dependent cases by following the line of arguments given in \cite{Fra25BOV} and \cite{FraKimtwocenterSZ}.  Recall that the system is governed by the Newtonian equation
 \begin{equation*} 
 \ddot{q}   = B(q)\dot{q} - \frac{q}{\lvert q \rvert^3} - \nabla E_t(q).
 \end{equation*}
whose solutions are critical points of the functional
\begin{equation*} 
\mathscr{A} \colon \mathfrak{L}\mathfrak{U} \to \R, \quad \mathscr{A}(q):=\frac{1}{2}\| \dot{q}\|^2 + \int_{S^1}q^* A - \int_0^1 V_t(q(t))dt.
\end{equation*}

Let $\mathfrak{Z}_0:=\Phi^{-1}(\mathfrak{U}_0)$ be an open subset of $\HH$ containing the origin. Note that $\mathfrak{Z}:=\mathfrak{Z}_0 \setminus \{ 0 \} = \Phi^{-1}(\mathfrak{U}).$ Since the zero set of every $z \in \mathfrak{L}\mathfrak{Z}$ is empty, the reparametrization map $t_z$ is now a diffeomorphism. Define
\begin{equation}\label{eq:SIGMA}
\Sigma \colon \mathfrak{L}\mathfrak{Z} \to \mathfrak{L}\mathfrak{U}, \quad \Sigma(z):=q_z,
\end{equation}
where $q_z$ is as in \eqref{eq:definitionofqz}.

We shall now  write the functional $\mathscr{A}$ in terms of $z,$ that is, pull it back via the map $\Sigma.$ We   assume that $z \in \mathfrak{L}\mathfrak{Z}$ satisfies
\begin{equation}\label{eq:conditoinofzinKS}
\left< z'(\tau), i z(\tau)\right> = - \left< i  z'(\tau), z(\tau) \right>=0, \quad \forall \tau \in S^1.
\end{equation}
We   compute 
\[
\dot{q}(t) = d\Phi( z(\tau_z(t)))z'(\tau_z(t))\dot{\tau}_z(t).
\]
The term $d\Phi(z)$ is explicitly computed  as 
\begin{equation*}
    d\Phi(z)v=\bar{v}iz+\bar{z}iv=-\overline{\bar{z}iv}+\bar{z}iv.
\end{equation*}
By conjugation we get 
\[
    d\Phi(z)^{\rm T}v=iz\bar{v}-izv,
\]
and for the special case $v\in\I\HH$ (i.e. $\bar{v}=-v$) this implies 
\begin{equation}\label{eq:dphitmatrix}
    d\Phi(z)^{\rm T}v=-2izv.
\end{equation}
Overall, the kinetic part becomes
\begin{align*}
   \ \int_0^1\frac{1}{2}\lvert \dot{q}(t)\rvert^2 dt &= \int_0^1 \frac{1}{2}\left< d\Phi( z(\tau_z(t)))z'(\tau_z(t)),  d\Phi(z(\tau_z(t))) z'(\tau_z(t))\right> \dot{\tau}_z(t)^2 dt \\
   &= \frac{\|z\|^2}{2} \int_0^1 \frac{\lvert \bar{z}'(\tau) i z(\tau) +\bar{z}(\tau) i z'(\tau) \rvert^2}{\lvert z(\tau)\rvert^2}d\tau \\
   &= 2 \|z\|^2 \|z'\|^2\\
   &=:\mathcal{K}(z),
\end{align*}
where the second identity follows from the change of variable $\tau =\tau_z(t),$ and the third identity follows from \eqref{eq:conditoinofzinKS}. 
The magnetic part is computed as:
\[
\mathcal{M}(z) = \int_{S^1}z^* \Phi^*A.
\]
Since
\[
 \int_0^1 E_t(q(t))dt = \frac{1}{\|z\|^2}\int_0^1 E_{t_z(\tau)}(\Phi(z(\tau)))\lvert z(\tau)\rvert^2 d\tau ,
 \]
the potential part is given by
\[
\mathcal{P}(z) = - \frac{1}{\| z\|^2} +\frac{1}{\|z\|^2}\int_0^1 E_{t_z(\tau)}(\Phi(z(\tau)))\lvert z(\tau)\rvert^2 d\tau .
\]
It follows that 
\begin{align*}
\Sigma^* \mathscr{A}(z) &= \mathcal{K}(z) + \mathcal{M}(z)- \mathcal{P}(z) \\
&= 2\|z\|^2\|z'\|^2 + \int_{S^1}z^* \Phi^*A+\frac{1}{\| z\|^2} -\frac{1}{\|z\|^2}\int_0^1 E_{t_z(\tau)}(\Phi(z(\tau)))\lvert z(\tau)\rvert^2 d\tau 
\end{align*}
for those $z \in \mathfrak{L}\mathfrak{Z}$ satisfying  \eqref{eq:conditoinofzinKS}. This functional naturally extends to a functional defined on  the open subset $\mathfrak{L}^*\mathfrak{Z}_0$, defined by
\begin{equation}\label{eq:L0Z0}
\mathfrak{L}^*\mathfrak{Z}_0:=\{ z \in \mathfrak{L}\mathfrak{Z}_0 \mid \|z\|\neq0\},
\end{equation}
via the same formula. We denote the extension by $\mathscr{B}\colon \mathfrak{L}^*\mathfrak{Z}_0 \to \R.$   Note that an element of $\mathfrak{L}^*\mathfrak{Z}_0$ is no longer required to   satisfy \eqref{eq:conditoinofzinKS}.
Further, note that  the regularized action functional $\mathscr{B}$ admits the $S^1$-action
\begin{equation}\label{eq:action}
\rho \colon S^1 \times \mathfrak{L}^*\mathfrak{Z}_0 \to \mathfrak{L}^* \mathfrak{Z}_0, \quad \quad \rho(\theta ,z) = e^{i\theta} z.
\end{equation}

Now we want to focus our attention on the critical points of the  regularized functional $\mathscr{B}.$ 
We first compute
\[
d\mathcal{K}(z)\xi = 4\|z'\|^2 \left<z,\xi\right>+4\|z\|^2\left<z',\xi'\right>= 4\|z'\|^2 \left<z,\xi\right>-4\|z\|^2\left<z'',\xi \right>.
\]
We  further  compute
\begin{align*}
    d\mathcal{M}(z)\xi &= \int_0^1\mathcal{L}_{\xi}(\Phi^*A)(z'(\tau))d\tau\\
    &= \int_0^1 ( \iota_{\xi}d\Phi^*A + d\iota_{\xi}\Phi^*A)(z'(\tau))d\tau\\
    &= \int_0^1 \Phi^*\sigma_B(z(\tau))(\xi(\tau), z'(\tau))d\tau\\
    &= \int_0^1 \sigma_B(\Phi(z(\tau)))( d\Phi(z(\tau))\xi(\tau), d\Phi(z(\tau))z'(\tau))d\tau\\
    &= \int_0^1  \left<      d\Phi(z(\tau)^{\rm T}  B(\Phi(z(\tau))) d\Phi(z(\tau))z'(\tau)),\xi(\tau)) \right>d\tau\\
    &=:    \left< \mathcal{N}(z)z', \xi\right> .
\end{align*}
Consider
\[
\mathcal{E}(z):=\frac{1}{\|z\|^2}\int_0^1 E_{t_z(\tau)}(\Phi(z(\tau)))\lvert z(\tau)\rvert^2 d\tau .
\]
To find its differential we first compute
\[
dt_z(\xi)(\tau) = \frac{2}{\|z\|^2}\int_0^\tau \left< z(s), \xi(s) \right>ds - \frac{2\left<z,\xi\right>}{\|z\|^4}\int_0^\tau \lvert z(s)\rvert^2 ds.
\]
Now we compute

\begin{align*}
    d\mathcal{E}(z)\xi&= - \frac{2\left<z,\xi\right>}{\|z\|^4}\int_0^1 E_{t_z(\tau)}(\Phi(z(\tau)))\lvert z(\tau)\rvert^2 d\tau\\
    &\;\;\;+ \frac{1}{\|z\|^2} \int_0^1\dot{E}_{t_z(\tau)}(\Phi(z(\tau)))dt_z(\xi)(\tau)\lvert z(\tau)\rvert^2d\tau\\
    &\;\;\; +\frac{1}{\|z\|^2}\int_0^1\left<  \nabla E_{t_z(\tau)}(\Phi(z(\tau))) ,  d\Phi(z(\tau))\xi(\tau)\right>\lvert z(\tau)\rvert^2d\tau \\
    &\;\;\; + \frac{2}{\|z\|^2}\int_0^1 E_{t_z(\tau)}(\Phi(z(\tau)))\left< z(\tau), \xi(\tau)\right>d\tau\\
    &=  - \frac{2\left<z,\xi\right>}{\|z\|^2}\mathcal{E}(z) +\frac{2}{\|z\|^4} \int_0^1\dot{E}_{t_z(\tau)}(\Phi(z(\tau)))\left(  \int_0^\tau \left< z(s), \xi(s) \right>ds\right)\lvert z(\tau)\rvert^2d\tau \\
    &\;\;\;- \frac{2\left<z,\xi\right>}{\|z\|^6}\int_0^1\dot{E}_{t_z(\tau)}(\Phi(z(\tau)))\left( \int_0^\tau \lvert z(s)\rvert^2 ds\right)\lvert z(\tau)\rvert^2d\tau\\
    &\;\;\; +\frac{1}{\|z\|^2}\int_0^1\lvert z(\tau)\rvert^2\left< d\Phi(z(\tau))^{\rm T}\cdot \nabla E_{t_z(\tau)}(\Phi(z(\tau))) ,  \xi(\tau)\right>d\tau \\
    &\;\;\; + \frac{2}{\|z\|^2}\int_0^1 E_{t_z(\tau)}(\Phi(z(\tau)))\left< z(\tau), \xi(\tau)\right>d\tau\\
        &=  - \frac{2\left<z,\xi\right>}{\|z\|^2}\mathcal{E}(z) + \frac{2}{\|z\|^4} \int_0^1 \left< z(s), \xi(s) \right> \left( \int_s^1 \dot{E}_{t_z(\tau)}(\Phi(z(\tau))) \lvert z(\tau)\rvert^2 d\tau \right)ds \\
    &\;\;\;- \frac{2\left<z,\xi\right>}{\|z\|^6}\int_0^1\dot{E}_{t_z(\tau)}(\Phi(z(\tau)))\left( \int_0^\tau \lvert z(s)\rvert^2 ds\right)\lvert z(\tau)\rvert^2d\tau\\
    &\;\;\; +\frac{2}{\|z\|^2}\int_0^1\lvert z(\tau)\rvert^2\left< -iz(\tau) \nabla E_{t_z(\tau)}(\Phi(z(\tau))) ,  \xi(\tau)\right>d\tau \\
    &\;\;\; + \frac{2}{\|z\|^2}\int_0^1 E_{t_z(\tau)}(\Phi(z(\tau)))\left< z(\tau), \xi(\tau)\right>d\tau 
\end{align*}
Therefore, if we write
\begin{align*}
\mathcal{E}^1(z)&=  \frac{1}{\|z\|^4} \int_0^1\dot{E}_{t_z(\tau)}(\Phi(z(\tau)))\left( \int_0^\tau \lvert z(s)\rvert^2 ds\right)\lvert z(\tau)\rvert^2d\tau   \\
\varepsilon_1(z)(\tau)&= \frac{1}{\|z\|^2}\left( \int_\tau^1 \dot{E}_{t_z(s)}(\Phi(z(s)))\lvert z(s)\rvert^2 ds\right)  z(\tau) \\
\varepsilon_2(z)(\tau)&=-iz(\tau) \nabla E_{t_z(\tau)}(\Phi(z(\tau)))\lvert z(\tau)\rvert^2 \\
    \varepsilon_3(z)(\tau)&= E_{t_z(\tau)}(\Phi(z(\tau)))z(\tau),
\end{align*}
then we have
\[
d\mathcal{P}(z)\xi = \frac{2}{\|z\|^4}\left< z , \xi \right> -\frac{2(\mathcal{E}(z) +\mathcal{E}^1(z))\left<z, \xi \right>}{\|z\|^2}+\frac{2\left< \varepsilon_1(z) + {\varepsilon}_2(z) +\varepsilon_3(z), \xi\right>}{\|z\|^2}
\]

Since every critical point $z \in \mathfrak{L}^* \mathfrak{Z}_0$ of $\mathscr{B}$ satisfies 
\[
d\mathscr{B}(z)\xi = d\mathcal{K}(z)\xi +d\mathcal{M}(z)\xi-d\mathcal{P}(z)\xi=0
\]
for every $\xi \in T_z \mathfrak{L}^*\mathfrak{Z}_0 = \mathfrak{L}\HH,$ it follows that $z$ satisfies the following second-order delay differential equation
\begin{equation}\label{eq:delaycriticalpoint}
\begin{aligned}
   z''(\tau)&= \left( \frac{\|z'\|^2}{\|z\|^2} +\frac{ \mathcal{E}(z) +\mathcal{E}^1(z)}{2\|z\|^4} - \frac{1}{2\|z\|^6}\right)z(\tau)+ \frac{1}{ 4\|z\|^2} \mathcal{N}(z)(\tau)z'(\tau)    \\
 &\;\;\;\;\;- \frac{\varepsilon_1(z)(\tau)+  \varepsilon_2(z)(\tau)+\varepsilon_3(z)(\tau)}{2\|z\|^4}
\end{aligned}
 \end{equation}

\begin{lemma} Suppose that $z \in \mathfrak{L}^*\mathfrak{Z}_0$ satisfies the delay equation \eqref{eq:delaycriticalpoint}. Then the zero set $\mathcal{Z}_z$ is finite.
\end{lemma}
\begin{proof} 
Fix a solution $z$ of \eqref{eq:delaycriticalpoint}. Then we can interpret the above delay equation as an ODE of the form
\[
z''(\tau) = C(\tau)z(\tau) + D(\tau)z'(\tau)
\]
for some $\tau$-dependent matrices $C$ and a  $D$. To see this note that according to \eqref{eq:delaycriticalpoint} all terms except $\epsilon_2$ only contribute a constant times the identity matrix to $C$. On the other hand $\epsilon_2$ can be interpreted as the matrix multiplication
\[
   \epsilon_2(z)(\tau)= M_{\epsilon_2}(z,\tau)z(\tau) = -iz(\tau) \nabla E_{t_z(\tau)}(\Phi(z(\tau)))\left< z(\tau),\ \cdot\ \right> (z(\tau)).
\]
Now let $\tau_*\in \mathcal{Z}_z$ be given and assume that $z'(\tau_*)=0$. Then by the Picard-Lindel\"of theorem the above linear ODE has a unique solution with the initial conditions:
\[
    z(\tau_*)=0 \quad \quad \text{ and } \quad \quad     z^\prime(\tau_*)=0
\]
It is easy to verify that this solution is $z\equiv0,$ contradicting the assumption $\|z \| \neq 0$ for all $ z\in\mathfrak{L}^* \mathfrak{Z}_0.$ Therefore, every zero of $z$ is transverse, from which the lemma follows.  \end{proof}

The previous lemma tells us that the reparametrization map $t_z$ as in \eqref{eq:tz} is well-defined for every critical point $z$ of $\mathscr{B},$ and it is a homeomorphism. Denote its continuous inverse by $\tau_z=t_z^{-1} \colon S^1 \to S^1.$ Define $q=q_z \colon S^1 \to \I\HH$ by the formula \eqref{eq:definitionofqz}, which is continuous on the whole circle and  smoothly differentiable on $S^1 \setminus t_z(\mathcal{Z}_z).$

\begin{theorem} \label{them:main11}
Let $z \in \mathfrak{L}^*\mathfrak{Z}_0$ be a critical point of the regularized functional $\mathscr{B}.$ Assume that it satisfies
\begin{equation}\label{eq:conditionofmainBOVtheorem}
\left< z'(\tau_0), i z(\tau_0)\right>=0 \quad \quad \text{ for some } \tau_0 \in S^1.
\end{equation}
Then  the map $q=q_z$ defined as in \eqref{eq:definitionofqz} is a generalized periodic solution of   the Newtonian equation \eqref{eq:secondorder}.  Moreover, the set of critical points of $\mathscr{B}$ modulo the $S^1$-action (see \eqref{eq:action}) is bijectively mapped onto the set of generalized solutions. 
\end{theorem}
\begin{proof} Assume that  $z \in \mathfrak{L}^*\mathfrak{Z}_0$ is a critical point of the functional $\mathscr{B},$ that is, it is a solution of the delay equation \eqref{eq:delaycriticalpoint}.  Let   $q=q_z \colon  S^1 \to \I \HH$ be the corresponding map given as in    \eqref{eq:definitionofqz}. The claim is straightforward in the case $\mathcal{Z}_z$ is empty. Therefore, in the following, without loss of generality, we may assume that $\#\mathcal{Z}_z=N $ for some positive integer $N \in \N. $

We then decompose the complement of the zero set into $N$ connected components
\[
S^1 \setminus \mathcal{Z}_z = \bigcup_{j=1}^N \mathcal{I}_j,
\]
where each component $\mathcal{I}_j = (\tau_j^-, \tau_j^+)$ is an open interval with the convention
\[
\tau_j^+ = \tau_{j+1}^-, \quad j \in \{1, \ldots, N-1\}\quad \text{ and } \quad \tau_N^+ = \tau_1^-.
\]
 The corresponding zero of $q=q_z$ is denoted by $t_j^\pm,$ that is,
 \[
 t_j^\pm := t_z(\tau_j^\pm), \quad 1 \leq j \leq N.
 \]
We set
\[
\mathcal{J}_j := t_z(\mathcal{I}_j)=(t_j^-, t_j^+)\subset S^1.
\]
Define
\[
\zeta(t):= z(\tau_z(t)),
\]
 then we have
 \begin{equation}\label{eq:qzetaandz}
 q(t) = \bar{\zeta}(t) i \zeta(t) \quad \quad \text{ and } \quad \quad \lvert q(t) \rvert = \lvert \zeta(t) \rvert^2 = \lvert z(\tau_z(t))\rvert^2. 
 \end{equation}
We compute
\begin{align*}
\frac{d}{d\tau}\left< z'(\tau), i z(\tau)\right> = & \; \left< z''(\tau) , i z(\tau)\right> \\
 = & \; -\frac{1}{2\|z\|^4}\left<z(\tau) \nabla E_{t_z(\tau)}(\Phi(z(\tau)))\lvert z(\tau)\rvert^2,z(\tau)\right>\\
    & \; + \frac{1}{4\|z\|^2 }\left<d\Phi(z(\tau))^{\rm T}  B(\Phi(z(\tau))) d\Phi(z(\tau))z'(\tau) ,iz(\tau)\right>\\
= & \; -\frac{1}{2\|z\|^4}\left< \nabla E_{t_z(\tau)}(\Phi(z(\tau)))\lvert z(\tau)\rvert^2,\lvert z(\tau)\rvert^2 \right>\\
    &\;  - \frac{1}{2\|z\|^2 }\left<  B(\Phi(z(\tau))) d\Phi(z(\tau))z'(\tau),\lvert z(\tau)\rvert^2 \right>\\
= & \;0,
\end{align*}
where we used (\ref{eq:delaycriticalpoint}) in the second identity. For the third identity we used the fact that for vectors restricted to $\mathbb{IH}$ the linear map $d\Phi(z(\tau))^{\rm T}$ coincides with the matrix representation of the quaternion $-2iz,$ see \eqref{eq:dphitmatrix}, and then we moved it to the other component of the vector product by transposing it.  
The last identity follows since $\nabla E_{t_z(\tau)}(\Phi(z(\tau)))$ and $B(\Phi(z(\tau))) d\Phi(z(\tau))z'(\tau)$ are both vectors in $\mathbb{IH}$, but $\lvert z(\tau)\rvert^2 $ has empty imaginary part. Hence, their vector products vanish.  
 Then condition \eqref{eq:conditionofmainBOVtheorem} implies that 
\begin{equation}\label{eq:globalconditionz}
\left< \dot{\zeta}(t), i \zeta(t)\right> = \dot{\tau}_z(t) \left< z'(\tau_z(t)), i z(\tau_z(t))\right>=0, \quad \forall t \in S^1.
\end{equation}
Restricting to the intervals $\mathcal{I}_j$ and $\mathcal{J}_j$ respectively, this  in particular implies 
\[
\dot{q}(t) = \dot{\bar{\zeta}}(t) i \zeta(t) + \bar{\zeta}(t) i \dot{\zeta}(t) = 2\bar{\zeta}(t) i \dot{\zeta}(t),
\]
and hence we obtain 
\[
\dot{\zeta}(t) =  - \frac{1}{2\lvert q(t) \rvert} i \zeta(t) \dot{q}(t). 
\]
We then find
\begin{align*}
z'(\tau_z(t)) =  {\dot{\zeta}(t)}t_z'(\tau) = \dot{\zeta}(t) \frac{\lvert z(\tau_z(t))\rvert^2}{\|z\|^2} = \dot{\zeta}(t) \frac{\lvert q(t)\rvert}{\|z\|^2} = - \frac{1}{2\|z\|^2} i \zeta(t) \dot{q}(t).
\end{align*}
We further differentiate:
\begin{align}
 \nonumber   z''(\tau_z(t))&= t_z'(\tau) \frac{d}{dt} z'(\tau_z(t))\\
\nonumber    &= \frac{ \lvert z(\tau_z(t))\rvert^2}{\|z\|^2} \frac{d}{dt}\left( - \frac{1}{2\|z\|^2} i \zeta(t) \dot{q}(t)\right)\\
 \nonumber   &= -\frac{1}{2\|z\|^4}\left( \frac{1}{2}\zeta(t) \dot{q}(t)^2 + \lvert q(t)\rvert i \zeta(t) \ddot{q}(t)\right)\\
 \label{eq:ztwoprime}   &= \frac{1}{2\|z\|^4}\left( \frac{1}{2}\lvert \dot{q}(t)\rvert^2 \zeta(t)- \lvert q(t)\rvert i \zeta(t) \ddot{q}(t)\right),
\end{align}
where we have used in the last identity the fact that $q \in \I \HH,$ so that $\dot{q}^2 = -\lvert \dot{q}\rvert^2.$

To plug \eqref{eq:delaycriticalpoint} into \eqref{eq:ztwoprime} we write each term in the right-hand side of \eqref{eq:delaycriticalpoint} in terms of $q$. We first introduce
 \[
 \mathfrak{C} \colon \mathfrak{L}\mathfrak{U}_0 \to \R, \quad \mathfrak{C}(q):=\int_0^1 E_t(q(t))dt
 \]
and
\[
\mathfrak{C}^1 \colon \mathfrak{L}\mathfrak{U}_0 \to \R, \quad \mathfrak{C}^1(q):= \int_0^1 t \dot{E}_t(q(t))dt.
\]
By the change of variable $\tau = \tau_z(t),$ we find
\[
\mathcal{E}(z) = \mathfrak{C}(q) \quad \quad \text{ and } \quad \quad \mathcal{E}^1(z)= \mathfrak{C}^1(q).
\]
Similarly, we find
\begin{align*}
 \varepsilon_1 (\tau_z(t)) &=    \left( \int_t^1 \dot{E}_s(q(s))ds \right) \zeta(t),\\
 \varepsilon_2(\tau_z(t)) &=    -i \zeta(t)  \nabla E_t(q(t))\lvert q(t)\rvert     ,\\
\varepsilon_3(\tau_z(t))&=   E_t(q(t)) \zeta(t). 
\end{align*}
Finally, we compute  
\begin{align*}
  \frac{1}{2}\|z\|^2 \mathcal{N}(z)(\tau)z'(\tau) = & \; \frac{1}{2}\|z\|^2d\Phi(z(\tau))^{\rm T}  B(\Phi(z(\tau))) d\Phi(z(\tau))z'(\tau)\\
  = & \; \frac{1}{2}\|z\|^2 (-2i \zeta(t) ) B(q(t)) d\Phi(\zeta(t))\dot{\zeta}(t)\cdot \frac{|q(t)|}{\|z\|^2}\\
  = &\; -  \lvert q(t)\rvert i\zeta(t)B(q(t)) \dot{q}(t)   
\end{align*}
Here we again used \eqref{eq:dphitmatrix} for the second identity.  
We then find
\begin{align*}
\frac{1}{2}\lvert \dot{q}(t)\rvert^2 \zeta(t) - \lvert q(t)\rvert i \zeta(t) \ddot{q}(t)  &=   -  \lvert q(t)\rvert i \zeta(t) B(q(t))\dot{q}(t) - \left( \int_t^1 \dot{E}_s(q(s))ds\right)\zeta(t)\\
&\;\;\;\;+C\zeta(t)+i\zeta(t) \nabla E_t(q(t)) \lvert q(t)\rvert - E_t(q(t)) \zeta(t)  
\end{align*}
where the constant $C$ is given by
\begin{equation}\label{eq:C}
    C = \frac{1}{2}\|\dot{q}\|^2 - \int_0^1 \frac{1}{\lvert q(t)\rvert}dt +\mathfrak{C}(q) +\mathfrak{C}^1(q).
\end{equation}
 Multiplying both sides by $\bar{\zeta}(t)i$ yields (recall that $q = \bar{\zeta} i \zeta$)
 \begin{align*}
\frac{1}{2}\lvert \dot{q}(t)\rvert^2 q(t) + \lvert q(t)\rvert^2   \ddot{q}(t) &=  \lvert q(t)\rvert^2 B(q(t))\dot{q}(t) - \left( \int_t^1 \dot{E}_s(q(s))ds\right)q(t)\\
&\;\;\;\;+Cq(t) + \nabla E_t(q(t)) \lvert q(t)\rvert^2 - E_t(q(t)) q(t).  
\end{align*}
We divide both sides by $\lvert q(t)\rvert^2$ and arrange this identity to obtain:
\begin{align*}
    \ddot{q}(t) - B(q(t)) \dot{q}(t) +\nabla E_t(q(t)) &=\left(C -\frac{1}{2}\lvert \dot{q}(t)\rvert^2       -   \int_t^1 \dot{E}_s(q(s))ds - E_t(q(t))\right)\frac{1}{   \bar{q}(t)   }
\end{align*}

We now argue as in the proof of \cite[Lemma 4.3]{FraKimtwocenterSZ}. We write
\[
\Phi(t):= C -\frac{1}{2}\lvert \dot{q}(t)\rvert^2  +\frac{1}{\lvert q(t)\rvert}     -   \int_t^1 \dot{E}_s(q(s))ds - E_t(q(t))
\]
so that the delay equation becomes
\begin{equation}\label{eeq:Phinitruou}
    \ddot{q}(t) - B(q(t)) \dot{q}(t) +\nabla E_t(q(t)) = \frac{1}{   \bar{q}(t)   } \Phi(t) - \frac{q(t)}{\lvert q(t)\rvert^3}
\end{equation}
We compute
\begin{align*}
    \dot{\Phi}(t) &= - \left< \dot{q}(t), \ddot{q}(t)\right> - \frac{ \left< q(t), \dot{q}(t)\right>}{\lvert q(t)\rvert^3} + \dot{E}_t(q(t)) - \dot{E}_t(q(t)) - \left< \nabla E_t(q(t)), \dot{q}(t)\right>\\
    &= - \left< \dot{q}(t), B(q(t))\dot{q}(t)   +\frac{1}{   \bar{q}(t)   }\Phi(t)  \right>   \\
        &= - \left< \dot{q}(t), \frac{1}{   \bar{q}(t)   }\Phi(t)  \right>   ,
\end{align*}
where we have used  \eqref{eeq:Phinitruou} in deriving the second identity and the skew-symmetry of $B$ in the third.  Since $\Phi(t)$ is real, it follows that it satisfies the ODE 
\[
\dot{\Phi}(t) = f(t)\Phi(t), \quad \quad t \in S^1 \setminus \mathcal{Z}_q,
\]
where
\[
f(t)= - \left< \dot{q}(t),  \frac{1}{   \bar{q}(t)   }  \right> 
\]
Writing this ODE in terms of $z$ and $\tau = \tau_z(t)$, we obtain
\[
\Phi'(\tau) = g(\tau) \Phi(\tau), \quad \quad \tau \in S^1 \setminus \mathcal{Z}_z,
\]
where
\begin{align*}
    g(\tau):= &\; f(\tau_z(\tau)) \frac{dt}{d\tau}\\
    =&  \;   - \left< \dot{q}(t),  \frac{1}{   \bar{q}(t)   }  \right>  \frac{dt}{d\tau}          \\
    =&  \;   - \left< \frac{2 i \|z\|^2z'(\tau)}{z(\tau)}, \frac{1}{ -{z}(\tau) i \bar{z}(\tau)} \right> \frac{\lvert z(\tau)\rvert^2}{\|z\|^2}          \\
    =&  \;     {\rm Re} \Bigg[  \frac{2  z'(\tau)}{z(\tau)} \Bigg] 
\end{align*}
We now fix arbitrary $\tau_0 \in \mathcal{I}_j$ and solve the ODE: we obtain
\[
\Phi(\tau) = \Phi(\tau_0) \frac{  \lvert z(\tau_0 )\rvert^2}{  \lvert z(\tau)\rvert^2} ,\quad \quad  \tau \in \mathcal{I}_j 
\]
from which we see that the function
\begin{align*}
    \Psi(\tau):= & \;   \Phi(\tau)  {\lvert z(\tau)\rvert^2}        \\
    =& \;    \lvert z(\tau)\rvert^2   \Bigg[ C - \int_0^1 \dot{E}_s (  \Phi(z(\tau_z(s)))       )  ds - E_{t_z(\tau)}( \Phi(z(\tau)  ))         \Bigg]  -2\| z\|^4 \lvert z'(\tau)\rvert^2 + 1                
\end{align*}
is constant on each connected component $\mathcal{I}_j.$ Note that $\Psi$ is continuous on the whole circle. It follows that there is a constant $K \in \R$ such that
\[
\Phi(\tau) = \frac{K}{\lvert z(\tau)\rvert^2}, \quad \quad \tau \in S^1 \setminus \mathcal{Z}_z,
\]
and thus, $\Phi$ does not change its sign. However, since
\[
\int_0^1 \Phi(t)dt =0,
\]
(using equation~\eqref{eq:C}) this implies that $\Phi$ vanishes identically and therefore $\Sigma(z)=q_z$ is a solution of the Newtonian equation \eqref{eq:secondorder}.

The last step of the proof is to show that the critical points of $\mathscr{B}$ modulo the $S^1$ action are in one to one correspondence to generalized solutions of \eqref{eq:secondorder}. For each $\theta \in S^1,$ we write $\rho_\theta=\rho(\theta, \cdot), $ where $\rho \colon S^1 \times \mathfrak{L}^*\mathfrak{Z}_0 \to\mathfrak{L}^*\mathfrak{Z}_0$ is given as in \eqref{eq:action}. 
For each $z \in \mathfrak{L}^*\mathfrak{Z}_0$ and $\xi \in T_z \mathfrak{L}^*\mathfrak{Z}_0 \equiv \mathfrak{L}\C,$ we have
 \[
 d\mathscr{B}(\rho_\theta (z))d\rho_\theta(z)\xi = d\mathscr{B}(z)\xi, 
 \]
from which we see that if $z$ is a critical point of $\mathscr{B},$ then $\rho_\theta(z)=e^{i \theta}z$ is a critical point as well. 

 Since we know that $\Sigma(z)$ is a solution of the ODE
 \eqref{eq:secondorder} we can use Lemma~\ref{S^1cover} to conclude that there is a one to one correspondence between the critical points of $\mathscr{B}$ modulo the $S^1$-action and the generalized solutions of \eqref{eq:secondorder}.  This completes the proof.  
\end{proof}

\begin{remark}
 In contrast to the planar cases \cite{Fra25BOV, FraKimtwocenterSZ} we cannot achieve that the regularized functional $\mathscr{B}$ is   Morse, due to the $S^1$-symmetry. The best we can get   for it is to be Morse-Bott.
\end{remark}

 Now we want to express this functional in the Hamiltonian formulation. For this we   apply the non-local Legendre transform which works analogously as the usual Legendre transform (see \cite[Section~7.1]{cieliebak2022a}).  

 First let us denote the non-local Lagrangian $ \mathscr{L} \colon \mathfrak{L}^*\mathfrak{Z}_0 \times \mathfrak{L}\HH \to \R $   by
\begin{align*}
    \mathscr{L}(z, v ) := & \; 2\|z\|^2\|v\|^2 - \int_0^1 2\left< izA(\bar{z}iz),v\right> +\frac{1}{\| z\|^2}\\ & \; -\frac{1}{\|z\|^2}\int_0^1 E_{t_z(\tau)}(\bar{z}iz)\lvert z(\tau)\rvert^2 d\tau. 
\end{align*} 
Note that the second term is just the one coming from $z^* \Phi^*A$.  Indeed, we compute
\[
z^*\Phi^* A = \left< A(\Phi(z )), d\Phi(z  )z' \right> d\tau =  \left< d\Phi(z)^{\rm T} A(\Phi(z )), z' \right> d\tau = -2\left< i zA(\bar{z}i z ), z'\right>d\tau,
\]
where we used the identity \eqref{eq:dphitmatrix} in the last step.
The corresponding action functional is given by
\[
\mathscr{L}(z,z')=\mathscr{B}(z), \quad z \in \mathfrak{L}^*\mathfrak{Z}_0.
\]
The next step is then to calculate the fiber derivative with respect to $\mathscr{L}$ which is the map
\begin{align*}
    F\mathscr{L} \colon  \mathfrak{L}^*\mathfrak{Z}_0 \times \mathfrak{L}\HH \to \mathfrak{L}^*\mathfrak{Z}_0 \times \mathfrak{L}\HH^*, \quad (z,v) \mapsto \left(z,
    \frac{d}{ds}\mathscr{L}(z,v+s\ \cdot\ )\right),
\end{align*}
where 
\[
    \frac{d}{ds}\mathscr{L}(z,v+s\ \cdot\ )\colon  u \mapsto \frac{d}{ds}\mathscr{L}(z,v+su)
\]
should be understood as an element in $\frak{L}\HH^*$. Computing this explicitly we get the new variable $w$:
\begin{align*}
     w:=   \frac{d}{ds}\mathscr{L}(z,v+s\ \cdot\ ) = \left<4\|z\|^2v-2izA(\bar{z}iz),\ \cdot\ \right>
\end{align*}
Since we are in $\HH$ we can simply identify the vector space with its dual vector space and write
\[
    w=4\|z\|^2v-2izA(\bar{z}iz).
\]
The inverse of the fiber translation is then simply given by
\[
    F\mathscr{L}^{-1}(z,w)= \left(z,\frac{w}{4\|z\|^2}+\frac{1}{2\|z\|^2}izA(\bar{z}iz)\right)
\]
With this we can now perform the non-local Legendre transform  as usual and obtain the non-local Hamiltonian:
\begin{align*}
    \mathcal{H}(z,w):= & \; \left<w,pr_2\circ F\mathscr{L}^{-1}(z,w)\right>-\mathscr{L}\left(F\mathscr{L}^{-1}(z,w)\right)\\
    = & \; \frac{\|w\|^2}{8\|z\|^2} + \frac{1}{2\|z\|^2}\left<w,izA(\bar{z}iz)\right>  +\frac{1}{2\|z\|^2}\left<izA(\bar{z}iz),izA(\bar{z}iz)\right> \\
    & \;-\frac{1}{\| z\|^2} +\frac{1}{\|z\|^2}\int_0^1 E_{t_z(\tau)}(\bar{z}iz)\lvert z(\tau)\rvert^2d\tau\\
    = & \; \frac{1}{8\|z\|^2}\|w+2izA(\bar{z}iz)\|^2-\frac{1}{\| z\|^2} +\frac{1}{\|z\|^2}\int_0^1 E_{t_z(\tau)}(\bar{z}iz)\lvert z(\tau)\rvert^2d\tau.
\end{align*}
The corresponding Hamiltonian action functional $\mathscr{A}_{\mathcal{H}} \colon \mathfrak{L}^*\mathfrak{Z}_0 \times \mathfrak{L} \HH \to \R$ is then defined as
\begin{align*}
    \mathscr{A}_{\mathcal{H}}(z,w):=   \left<w, z^\prime \right> - \mathcal{H}(z,w)= \int_{S^1} (z,w)^*\lambda - \mathcal{H}(z,w),
    \end{align*}
    where $\lambda=\sum\limits_{i=1}^4 w_idz_i$ denotes the canonical one-form. Computing the differential of this functional we have
    \begin{align*}
        d\mathscr{A}_{\mathcal{H}}(z,w)(\varphi_z,\varphi_w)= & \left<\varphi_w,z^\prime\right> + \left<w,\varphi_z^\prime\right> - \left<\nabla_z\mathcal{H}(z,w),\varphi_z\right> - \left<\nabla_w\mathcal{H}(z,w),\varphi_w\right>\\
        = & \left<z^\prime-\nabla_w\mathcal{H}(z,w),\varphi_w\right> - \left<w^\prime+\nabla_z\mathcal{H}(z,w),\varphi_z\right>.
    \end{align*}
    Consequently, the critical points of the Hamiltonian action functional $\mathscr{A}_{\mathcal{H}}$ are the $1$-periodic solutions of the  Hamiltonian equations of $\mathcal{H},$ which   are given by
    
\begin{equation}\label{eq:nonlocalHameq}
    \begin{aligned}
        z' = \;\; \nabla_{w} \mathcal{H} \; = &  \; \frac{1}{4\|z\|^2}\left(w+2izA(\bar{z}iz)\right),    \\
        w' = -  \nabla_z \mathcal{H} =  & \; -\frac{1}{4\|z\|^4} \|w+2izA(\bar{z}iz)\|^2z- \frac{2z}{\|z\|^4} \\ & \;- \frac{1}{4\|z\|^2}\left(2iA(\bar{z}iz)-\mathcal{N}(z)\right)\left(w+2izA(\bar{z}iz)\right)\\ &\; +\frac{2(\mathcal{E}(z) +\mathcal{E}^1(z))z}{\|z\|^2}-\frac{2( \varepsilon_1(z) + {\varepsilon}_2(z) +\varepsilon_3(z))}{\|z\|^2}     
    \end{aligned}
\end{equation}
Note that the critical points of $\mathscr{A}_{\mathcal{H}}$ are in one to one correspondence to the critical points of $\mathscr{B}$ via the fiber translation $F\mathscr{L}$ (see \cite[Proposition~7.1]{cieliebak2022a} for more details).

\begin{remark}
     The non-local Lagrangian and Hamiltonian action functionals $\mathscr{B}$ and $\mathscr{A}_{\mathcal{H}}$ satisfies 
    \[
    \mathscr{A}_{\mathcal{H}}(z,w) = \mathscr{B}(z) - \frac{\| w - 4\|z\|^2 z'\|^2}{8\|z\|^2}  - \frac{\left< w + izA(\bar{z}i z) - 4\|z\|^2z', i z A(\bar{z}i z)\right>}{2\|z\|^2}.
    \]
    In particular, they    coincide on critical point, that is, if $(z, w)$ satisfies \eqref{eq:nonlocalHameq}, then we have
    \[
    \mathscr{A}_{\mathcal{H}}(z, w) = \mathscr{B}(z).
    \]
\end{remark}

\bibliographystyle{abbrv}
\bibliography{mybibfile}

\end{document}